\documentclass{amsart}

\newtheorem{theorem}{Theorem}[section]

\theoremstyle{definition}

\newtheorem{example}[theorem]{Example}

\theoremstyle{remark}

\begin{document}

\title[Numerical Landen]{An 
iterative method for numerical integration of rational functions}

\author{Dante Manna}
\address{Department of Mathematics and Statistics, Dalhousie University, 
Halifax, Nova Scotia, Canada B3H 3J5}
\email{dmanna@mathstat.dal.ca}

\author{Victor H. Moll}
\address{Department of Mathematics,
Tulane University, New Orleans, LA 70118}
\email{vhm@math.tulane.edu}

\subjclass{Primary 33}

\date{\today}


\begin{abstract}
We describe a new method for numerical integration of rational  functions
on the real line. Given a rational integrand, we provide a new rational 
function preserving its integral on the line. The coefficients 
of the new function 
are explicit polynomials in the original ones. These transformations 
depend on the degree of the input and the desired order of the method.
Both parameters are arbitrary. The 
formulas can be precomputed. Iteration yields an approximation of the 
desired integral, with $m$-th order convergence. Examples  illustrating the 
automatic 
generation of these formulas and a comparison with standard numerical schemes
are also presented.
\end{abstract}

\maketitle

\newcommand{\nn}{\nonumber}
\newcommand{\ba}{\begin{eqnarray}}
\newcommand{\ea}{\end{eqnarray}}
\newcommand{\ift}{\int_{0}^{\infty}}
\newcommand{\ifft}{\int_{- \infty}^{\infty}}
\newcommand{\eps}{\epsilon}
\newcommand{\no}{\noindent}
\newcommand{\realpart}{\mathop{\rm Re}\nolimits}
\newcommand{\imagpart}{\mathop{\rm Im}\nolimits}

\newtheorem{Definition}{\bf Definition}[section]
\newtheorem{Thm}[Definition]{\bf Theorem}
\newtheorem{Example}[Definition]{\bf Example}
\newtheorem{Lem}[Definition]{\bf Lemma}
\newtheorem{Note}[Definition]{\bf Note}
\newtheorem{Cor}[Definition]{\bf Corollary}
\newtheorem{Prop}[Definition]{\bf Proposition}
\newtheorem{Problem}[Definition]{\bf Problem}
\numberwithin{equation}{section}

\section{Introduction} \label{sec-intro}
\setcounter{equation}{0}

The problem of numerical integration of a function over the real line 
is described in Numerical Analysis texts such 
as  \cite{davisrabi1, ueberhuber1}. The standard algorithms for the 
numerical integration of 
\begin{equation}
I := \int_{- \infty}^{\infty} F(x) \, dx
\end{equation}
\noindent 
start by dealing with the unboundedness of the domain of integration. This is 
usually resolved in two ways: the first one considers the problem on a 
finite interval 
\begin{equation}
I_{L} := \int_{-L}^{L} F(x) \, dx 
\end{equation}
\noindent 
followed by a convergence study as $L \to \infty$. The alternative is to 
transform the real line to a bounded interval. For example, the map 
$t = x/(1+x)$ maps $[0, \infty)$ to $[0,1]$ and then 
\begin{eqnarray}
\int_{-\infty}^{\infty} F(x) \, dx  & =  & 
\ift \left( F(x) + F(-x) \right) \, dx  \label{newint} \\
& = & \int_{0}^{1} \left[ F \left( \frac{t}{1-t} \right) + 
F \left(  \frac{t}{t-1} \right) \right] \frac{dt}{(1-t)^{2}}. \nonumber 
\end{eqnarray}
\noindent
The unboundedness of the original interval of integration is now reflected 
in the  (possible) singularity of the new integrand at the boundary $t=1$. 
Observe that if the original integrand is a rational function, then so is 
the new one in (\ref{newint}). 

In this paper we present a new numerical method for the integration of 
rational functions on $\mathbb{R}$. It is different in spirit to the 
standard ones: the numerical approximation to the integral is obtained from
a {\em recurrence} acting on the coefficients of the integrand. In
particular, the integrand is never evaluated. We illustrate the 
comparison with the standard methods; a more systematic study will be 
presented elsewhere. 

The method presented here is based on a transformation of the coefficients of 
a rational function, that preserves its integral. This is the so-called 
\emph{rational Landen transformation}. The original problem discussed by 
Landen, Gauss, Legendre and others deals with the elliptic integral 
\begin{equation}
G(a,b) = \int_{0}^{\pi/2} \frac{d \theta}{\sqrt{a^{2} \cos^{2} \theta  + 
b^{2} \sin^{2} \theta } }.
\end{equation}
\noindent 
Gauss \cite{gauss1} made the {\em numerical observation} that the function 
$G(a,b)$ was invariant under the transformation 
\begin{equation}
(a,b ) \mapsto \left( \frac{a+b}{2},  \sqrt{ab} \right). 
\label{agm}
\end{equation}
\noindent
The iteration of (\ref{agm}) leads to a 
numerical evaluation of the elliptic integral, or, as it has been 
explained  by J. and P. Borwein in \cite{borwein1}, to efficient methods
for the numerical evaluation of $\pi$. The rational 
analogue of this transformation was developed in \cite{boros1, manna-moll2}
and here we show how to use it as a numerical method to evaluate rational 
integrals. The reader will find in \cite{manna-moll3} a survey of the diverse
aspects related to these transformations.  

Section 2 discusses the basic structure of the algorithm. Section 3 introduces 
a family of polynomials that play an important role in the development of the
formulas given in Section 4. Finally, Section 5 contains some examples. The 
first illustrates the steps for a method of order $2$ acting on a rational 
function of degree $6$. The next two examples illustrate the accuracy of the 
method and its comparison to the trapezoidal rule. A systematic study of the 
cost involved in this algorithm will be presented elsewhere. 

\section{The Landen transformation and algorithm}
\setcounter{equation}{0}

We present a general description of an iterative algorithm for the evaluation of
\begin{equation}
I := \int_{-\infty}^{\infty} F(x) \, dx.
\end{equation}
\noindent 
Here $F$ is a rational function given as 
\begin{equation}
F(x) = \frac{B(x)}{A(x)},
\label{F-def}
\end{equation}
\noindent
with
\begin{equation}
A(x) = \sum_{k=0}^{p} a_{k}x^{p-k} \text{ and } 
B(x) = \sum_{k=0}^{p-2} b_{k}x^{p-2-k}, 
\end{equation}  
where $a_{p} \neq 0$.  The general 
construction treats the coefficients $a_i$ and $b_j$ as indeterminates.  
Naturally, for specific integrands the parameters $a_i$ and $b_j$ are 
real numbers and the maps described in this section are defined on 
parts of ${\mathbb{R}}^{2p}$ where the integrals are convergent. \\

The set
\begin{equation}
{\mathcal{R}}^{2p}: 
\{ (a_{0}, a_{1}, \ldots, a_{p}, b_{0}, \ldots, b_{p-2}): \text{ 
such that } F \text{ in } (\ref{F-def}) \text{ has finite integral } 
\}, 
\nonumber
\end{equation}
\noindent
will be used to represent the rational function $F$ in terms of its 
coefficients. It will be referred as the
\emph{coefficient space}.

\medskip

In Section  \ref{sec-algo} we describe the 
construction, for each integer $m \geq 2$, of a rational function 
$F_{1,m}$ that
satisfies 
\begin{equation}
\int_{-\infty}^{\infty} F_{1,m}(x) \, dx = 
\int_{-\infty}^{\infty} F(x) \, dx. \label{invariance}
\end{equation}
The function $F_{1,m}(x)$ has the same degree as the original $F$
and the new coefficients are \emph{polynomials} in the  old coefficients 
$a_0, \ldots,a_p, b_0,\ldots,b_{p-2}$. 
Naturally, this produces a map on ${\mathcal{R}}^{2p}$ that we denote by 
$\mathfrak{L}_{m,p}:  \mathcal{R}^{2p} \rightarrow \mathcal{R}^{2p}$, called 
the \emph{rational Landen 
transformation of order $m$ and degree $p$}. We denote by 
$\mathfrak{L}_{m,p}^n$ its $n$-fold composition. Now introduce a new operator
by $\phi : \mathcal{R}^{2p} \rightarrow \mathbb{R}$ as the map which 
takes a vector of length $2p$ corresponding to a rational function 
in $\mathcal{R}$ of degree $p$ and returns the value of the rational function 
at $x=0$. In terms of the vector of coefficients, this is the ratio 
of the last entry  over its 
$(p+1)-$th one. The composition $\phi \circ \mathfrak{L}_{m,p}^{n}$ will 
be denoted by $\phi_{m,p}^{n}$.  The motivation behind $\phi$ is 
the following: if 
\begin{equation}
F(x) = \frac{b_{0}+ b_{1}x + \cdots + b_{p-2}x^{p-2} }
{a_{0}+ a_{1}x + \cdots + a_{p}x^{p} }, \label{F}
\end{equation}
\noindent
then the function obtained by iterating the  map 
$\mathfrak{L}_{m,p}$ $n$ times is written as 
\begin{equation}
F_{n,m}(x) = \frac{b_{0,n}+ b_{1,n}x + \cdots + b_{p-2,n}x^{p-2} }
{a_{0,n}+ a_{1,n}x + \cdots + a_{p,n}x^{p} }. 
\end{equation}
\noindent
The value of $F_{n,m}$ at $x=0$ gives a sequence of real numbers that 
converges to $1/\pi$ times the integral of $F$ in (\ref{F}); see 
\cite{manna-moll2} for details. The update on this sequence comes from 
the coefficients of $F_{n,m}$. These are obtained by 
applying $\mathfrak{L}_{m,p}$
to those of $F_{n-1,m}$. 
Finally, define $\vec{\alpha} := (a_0,...,a_p,b_0,...,b_{p-2})$.

The main result of \cite{manna-moll2} is that the Landen transformation 
satisfies 
\begin{equation} 
\phi^n_{m,p}(\vec{\alpha}) \rightarrow \frac{I}{\pi} 
\end{equation} 
as $n \to \infty$ if $I < \infty$.  Furthermore, the convergence is of
\emph{order} $m$, that is, 
\begin{equation} 
\left| \phi^{n+1}_{m,p}(\vec{\alpha}) - \frac{I}{\pi}\right| 
\leq C \left|\phi^n_{m,p}(\vec{\alpha})- \frac{I}{\pi} \right|^m.
\label{2point8}
\end{equation} 
This convergence result appears in \cite{hubbard1} for the case $m=2$ 
and the general case can be established along the same lines. See 
\cite{manna1} for details. 

\medskip

Adapting the rational Landen transformations into a numerical method for 
calculating $I$ involves a process with two parts.  The first one 
is a symbolic calculation of the explicit algebraic formulae for the rational 
Landen transformation.  The steps in this calculation are described in 
detail in Section \ref{sec-algo}. The second one is the iteration of these 
formulas. 

\bigskip

\noindent
{\bf Algorithm 1}

\medskip

\noindent
{\bf Input}: 

\noindent
1) An integer $p$: the degree of the denominator $A$. 

\noindent
2) An integer $m \geq 2$:  the order of the transformation. 

\medskip

\noindent
{\bf Output}:

The explicit formula for the transformation 
$\mathfrak{L}_{m,p}: \mathcal{R}^{2p} \rightarrow \mathcal{R}^{2p}$, as 
polynomials in the $a_i$ and $b_j$.

\medskip

\begin{Note}
Observe that, given $m$ and $p$, the map $\mathfrak{L}_{m,p}$ can be 
\emph{precomputed} and the result can be stored for its use in the 
second algorithm. Therefore, the first 
algorithm carries a \emph{one-time cost} and is not figured into the 
time of the method. This 
precomputation will be assumed in the discussion of the second 
algorithm. 
\end{Note}

\medskip

\noindent
{\bf Algorithm 2}

\medskip

\noindent
{\bf Input}: 

\noindent
1) A vector $\vec{\alpha}$ representing the coefficients of the rational 
integrand $F$.

\noindent
2) An integer $m \geq 2$:  the order of convergence.

\noindent
3) An integer $n \in \mathbb{N}$: the number of iterations of the Landen map. 

\medskip

\noindent
{\bf Output}:
 The expression $\phi^n_{m,p}(\vec{\alpha})$ that approximates $I/\pi$.

\section{The evaluation of the polynomials $P_{m}$ and $Q_{m}$} 
\label{sec-prepa}
\setcounter{equation}{0}

The algorithm described in Section \ref{sec-algo} employs the polynomials 
\begin{equation}
P_{m}(x) := \sum_{j=0}^{\lfloor{ m/2 \rfloor}} (-1)^{j} \binom{m}{2j}x^{m-2j} 
\label{Pm}
\end{equation}
\noindent
and
\begin{equation}
Q_{m}(x) := \sum_{j=0}^{\lfloor{(m-1)/2 \rfloor}} (-1)^{j} \binom{m}{2j+1}
x^{m-(2j+1)}. 
\label{Qm}
\end{equation}
\noindent
The integration algorithm discussed here is based on the fact that the 
rational function 
\begin{equation}
R_{m}(x) := \frac{P_{m}(x)}{Q_{m}(x)}
\end{equation}
\noindent
satisfies 
\begin{equation}
\cot(m \theta) = R_{m}(\cot \theta).
\end{equation}
\noindent
For instance, for $m=2$, we have
\begin{equation}
P_{2}(x) = x^{2}-1, \text{ and } Q_{2}(x) = 2x. 
\end{equation}

\begin{Note}
This trigonometric property is instrumental in the proof that the integral of 
$F$ is the same as that of $F_{1}$. See \cite{manna-moll2} for details.
\end{Note}

\section{The algorithm} \label{sec-algo}
\setcounter{equation}{0}

In this section we describe each of the steps in the first algorithm. This 
algorithm has been implemented in Mathematica 6.0.  \\

\noindent
{\bf Step 1}. Construct the polynomial
\begin{equation}
H(x) = \sum_{i=0}^{p} h_{i}x^{p-i} 
\end{equation}
\noindent
defined by 
\begin{equation}
H(x) := \text{Res}_{z} \left( A(z), P_{m}(z) - x Q_{m}(z) \right), 
\label{res1}
\end{equation}
\noindent
where $\text{Res}_{z}$ denotes the \emph{resultant} in the variable $z$. The 
degrees of the 
polynomials involved in (\ref{res1}) are $p= \text{deg}A$ and 
$m = \text{deg}(P_{m}(z)-xQ_{m}(z))$, respectively. \\

\noindent
{\em The resultant}. Given two polynomials $\alpha(t)$ and 
$\beta(t)$, the resultant
of $\alpha$ and $\beta$ is defined by
\begin{equation}
\text{Res}(\alpha,\beta) := \prod_{i=1}^{r} \prod_{j=1}^{s} (y_{j} - x_{i}),
\end{equation}
\noindent
where $x_{i}$ are the roots of $\alpha(t)=0$ and $y_{j}$ are 
the roots of $\beta(t)=0$. 

The resultant of two polynomials can be computed as the determinant of the 
\emph{Sylvester matrix} formed by their coefficients; see 
\cite{lang-algebra}. For instance, if 
\begin{equation}
\alpha(t) = a_{0} + a_{1}t + a_{2}t^{2} + a_{3}t^{3}  \text{ and }
\beta(t) = b_{0} + b_{1}t + b_{2}t^{2},
\nonumber
\end{equation}
\noindent
then the Sylvester matrix is defined by 
\begin{equation}
S_{3,2} := \begin{pmatrix} 
a_{3} & a_{2} & a_{1} & a_{0} & 0 & 0 \\
  0   & a_{3} & a_{2} & a_{1} & a_{0} & 0 \\
  0   & 0  & a_{3} & a_{2} & a_{1} & a_{0} \\
  b_{2} & b_{1} & b_{0} & 0 & 0 & 0 \\
  0 & b_{2} & b_{1} & b_{0} & 0 & 0  \\
  0 & 0  & b_{2} & b_{1} & b_{0} & 0   \\
  0 & 0 & 0 & b_{2} & b_{1} & b_{0}
\end{pmatrix}
\label{example1}
\end{equation}
\noindent
and it is a square matrix of 
size $\text{deg}(\alpha) + \text{deg}(\beta) + 2 = 7$. 
The resultant of $\alpha(t)$ and $\beta(t)$ in (\ref{example1}) is 
\begin{eqnarray}
\text{Res}(\alpha, \beta) & = & a_{3}^{2}b_{0}^{3} -a_{2}a_{3}b_{0}^{2}b_{1} 
+ a_{1}a_{3}b_{0}b_{1}^{2}-a_{0}a_{3}b_{1}^{3} + 
a_{2}^{2}b_{0}^{2}b_{2} \nonumber \\
& & -2a_{1}a_{3}b_{0}^{2}b_{2} - a_{1}a_{2}b_{0}b_{1}q_{2} + 
3a_{0}a_{3}b_{0}b_{1}b_{2}+a_{0}a_{2}b_{1}^{2}b_{2} \nonumber \\
& & +a_{1}^{2}b_{0}b_{2}^{2} - 2a_{0}a_{2}b_{0}b_{2}^{2} - 
a_{0}a_{1}b_{1}b_{2}^{2} + a_{0}^{2}b_{2}^{3}. \nonumber 
\end{eqnarray}

In general, the resultant of two polynomials is a polynomial in their 
coefficients. 

\begin{Note}
The polynomial $H$ has the same degree as $A$, the denominator of the integrand 
$R(x)$. It will become the denominator of the new rational function. 
Its coefficients $h_{i}$ are polynomials in those of 
$A$. The calculation of $H$ can be obtained by evaluating the determinant
of a square matrix of dimension 
$\text{deg}(A) + \text{deg}(P_{m}) +2 = p+m+2$. For instance, for 
$m=2$ and 
\begin{equation}
A(x) = a_{0}x^{4}+a_{1}x^{3} + a_{2}x^{2}+a_{1}x + a_{4}, 
\end{equation}
\noindent
we obtain
\begin{eqnarray}
H(x) & = & 16a_{0}a_{4}x^{4} + 8(a_{1}a_{4}-a_{0}a_{3})x^{3} + 
4(a_{0}a_{2}-a_{1}a_{3} 
 +  4a_{0}a_{4}+a_{2}a_{4})x^{2} \nonumber \\
& + & 2(-a_{0}a_{1}+a_{1}a_{2}-3a_{0}a_{3}-a_{2}a_{3}+3a_{1}a_{4}+a_{3}a_{4})x 
\nonumber \\
& + & (a_{0}-a_{1}+a_{2}-a_{3}+a_{4})(a_{0}+a_{1}+a_{2}+a_{3}+a_{4}). 
\nonumber
\end{eqnarray}
\end{Note}

\medskip

\noindent
{\bf Step 2}. Form the polynomial 
\begin{equation}
E(x) = \sum_{i=0}^{p} h_{i} \left( P_{m}(x) \right)^{p-i} \, 
\left( Q_{m}(x) \right)^{i}.
\label{def-E}
\end{equation}

\medskip

\noindent
{\bf Step 3}. The polynomial $E(x)$ formed in Step 2 is a multiple of the 
denominator $A(x)$. Compute the quotient 
\begin{equation}
Z(x) = \frac{E(x)}{A(x)} 
\end{equation}
\noindent
and write it as
\begin{equation}
Z(x) = \sum_{k=0}^{r} z_{k}x^{r-k}, \text{ with } r = p(m-1).
\end{equation}

\noindent
{\bf Step 4}. Compute the product
\begin{equation}
C(x) = B(x) Z(x) 
\end{equation}
\noindent
and write it as 
\begin{equation}
C(x) = \sum_{k=0}^{s} c_{k}x^{s-k}, \text{ with } s = mp-2. 
\end{equation}

\medskip

\noindent
{\bf Step 5}. Form the expression 
\begin{equation}
T_{x}(a,b) := \sum_{j=0}^{x} (-1)^{a-x+j} \binom{a}{x-j} \binom{b}{j}, 
\end{equation}
\noindent
for $a, \, b, \, x \in \mathbb{N}$. 

\medskip

\noindent
{\bf Step 6}. Define the expressions
\begin{eqnarray}
M_{1}(j, \alpha, \beta, \gamma, m,p)  & :=  & 
(-1)^{j+ \alpha - \beta} c_{2j} \frac{2^{2(\alpha- \beta)} \alpha}
{2 \alpha - \beta} \binom{2 \alpha - \beta}{\beta} 
\binom{\nu - \alpha -1 + \beta}{\gamma}  \times \nonumber \\
 & \times & 
\left[ T_{\lambda + \alpha m }(2j, s - 2j) + 
T_{\lambda - \alpha m }(2j, s - 2j)  \right], \nonumber
\end{eqnarray}
\noindent
and 
\begin{eqnarray}
M_{2}(j, \alpha, \beta, \gamma, m,p)  & :=  & 
(-1)^{j+ \beta} c_{2j+1} 2^{2 \beta +1}
\binom{\alpha + \beta}{2 \beta +1} 
\binom{\nu - 2 - \beta}{\gamma}  \times \nonumber \\
 & \times & 
\left[ T_{\lambda + \alpha m }(2j+1, s - 2j-1) - 
T_{\lambda - \alpha m }(2j+1, s - 2j-1)  \right], \nonumber
\end{eqnarray}
\noindent
with $\nu := p/2$ and $\lambda := (mp-2)/2$.  \\

\noindent
{\bf Step 7}. Define 
\begin{eqnarray}
J(x) & := & \frac{1}{2^{s}} 
\sum_{\gamma =0}^{\nu-1} \left( \binom{\nu-1}{\gamma} 
\sum_{j=0}^{\lambda} (-1)^{j} c_{2j} T_{\lambda}(2j,s-2j) \right)x^{2 \gamma} + 
\nonumber \\
 & + & \frac{1}{2^{s}} 
\sum_{\gamma=0}^{\nu-2} \left( \sum_{j=0}^{\lambda} 
\sum_{\alpha=1}^{\nu-1-\gamma} \sum_{\beta=0}^{\alpha} 
M_{1}(j, \alpha, \beta, \gamma, m,p) \right) x^{2 \gamma} \nonumber \\
 & + & \frac{1}{2^{s}} 
\sum_{\gamma=1}^{\nu-1} \left( \sum_{j=0}^{\lambda} 
\sum_{\alpha=\nu-\gamma}^{\nu-1} \sum_{\beta=\alpha-\nu+\gamma+1}^{\alpha} 
M_{1}(j, \alpha, \beta,\gamma,  m,p) \right) x^{2 \gamma} \nonumber \\
 & + & \frac{1}{2^{s}} 
\sum_{\gamma=0}^{\nu-2} \left( \sum_{j=0}^{\lambda-1} 
\sum_{\alpha=1}^{\nu-1-\gamma} \sum_{\beta=0}^{\alpha-1} 
M_{2}(j, \alpha, \beta,\gamma,  m,p) \right) x^{2 \gamma+1} \nonumber \\
 & + & \frac{1}{2^{s}} 
\sum_{\gamma=1}^{\nu-2} \left( \sum_{j=0}^{\lambda-1} 
\sum_{\alpha=\nu-\gamma}^{\nu-1} \sum_{\beta=0}^{\alpha-1} 
M_{2}(j, \alpha, \beta,\gamma,  m,p) \right) x^{2 \gamma+1}. \nonumber 
\end{eqnarray}

\medskip
\noindent
{\bf Step 8}. The new rational function is defined by
\begin{equation}
F_{1,m}(x) := \frac{J(x)}{H(x)}. 
\end{equation}
\noindent
It satisfies (\ref{invariance}). \\

The reader is referred to \cite{manna-moll2} for the proofs of the formulas 
describing this algorithm. 

\section{Examples} \label{sec-examples}
\setcounter{equation}{0}

In this section we give examples that illustrate the rational Landen 
transformations. 

\begin{example}
We provide a step by step construction of the rational 
Landen transformation of order $m=2$ for the function
\begin{equation}
F(x) = \frac{b_{0}x^{4} + b_{1}x^{3} + b_{2}x^{2} + b_{3}x + b_{4} }
{a_{0}x^{6} + a_{1}x^{5} + a_{2}x^{4} + a_{3}x^{3} + a_{4}x^{2} + a_{5}x + a_{6} }. 
\end{equation}
\noindent
The goal is to produce a function $F_{1,6}(x)$ with the same integral as 
$F(x)$.

\medskip

\begin{Note}
The special case $b_{1} = b_{3} = a_{1} = a_{3} = a_{5} = 0$ was the first 
example of this new type of transformation. It appears in \cite{boros2}.
\end{Note}

\begin{Note}
The choice of $m=2$ requires the evaluation of the
polynomials $P_{2}(x) = x^{2}-1$ and $Q_{2}(x) = 2x$. These  are computed 
directly from (\ref{Pm}) and (\ref{Qm}).  
\end{Note}

\medskip

\noindent
{\bf Step 1} computes the polynomial $H(x)$ from (\ref{res1}). The Mathematica 
command \emph{Resultant} yields the expressions
\begin{eqnarray}
h_{0} & = & 64a_{0}a_{6} \nonumber  \\
h_{1} & = & -32(a_{0}a_{5}-a_{1}a_{6}) \nonumber  \\
h_{2} & = & 16(a_{0}a_{4}-a_{1}a_{5}+6a_{0}a_{6} + a_{2}a_{6}) \nonumber  \\
h_{3} & = & -8(a_{0}a_{3} -a_{1}a_{4}+5a_{0}a_{5}+a_{2}a_{5}-5a_{1}a_{6}-a_{3}a_{6}) \nonumber  \\
h_{4} & = & 4(a_{0}a_{2}-a_{1}a_{3}+4a_{0}a_{4}+a_{2}a_{4}-4a_{1}a_{5}-a_{3}a_{5}+9a_{0}a_{6}+4a_{2}a_{6}+a_{4}a_{6}) \nonumber \\
h_{5} & = & -2(a_{0}a_{1}-a_{1}a_{2}+3a_{0}a_{3}+a_{2}a_{3}-3a_{1}a_{4}-
a_{3}a_{4}+5a_{0}a_{5}) \nonumber \\
 & & - 2(3a_{2}a_{5}+a_{4}a_{5}-5a_{1}a_{6}-3a_{3}a_{6} - 
a_{5}a_{6}) \nonumber \\
h_{6} & = & (a_{0}-a_{1}+a_{2}-a_{3}+a_{4}-a_{5}+a_{6})
(a_{0}+a_{1}+a_{2}+a_{3}+a_{4}+a_{5}+a_{6}). \nonumber 
\end{eqnarray}
\noindent
The polynomial $H(x)$ is the denominator of the new rational function obtained 
as a product of the Landen transformation. 

\medskip

\noindent
{\bf Step 2} computes the polynomial $E(x)$ from (\ref{def-E}). In this 
example, this is a polynomial of degree $12 (= mp)$, that we write as 
\begin{equation}
E(x) = \sum_{i=0}^{12} e_{i}x^{12-i}.
\end{equation}
\noindent
The  symbolic expansion of (\ref{def-E}) produces 
\begin{eqnarray}
e_{0} & = & e_{12} = 64a_{0}a_{6} \nonumber \\
e_{1} & = & - e_{11} = 64(a_{0}a_{5}-a_{1}a_{6}) \nonumber \\
e_{2} & = & e_{10} = 64(a_{0}a_{4}-a_{1}a_{5}+a_{2}a_{6}) \nonumber \\
e_{3} & = & -e_{9} = 64(a_{0}a_{3}-a_{1}a_{4}+a_{2}a_{5}-a_{3}a_{6}) \nonumber \\
e_{4} & = & e_{8} = 64(a_{0}a_{2}-a_{1}a_{3} + a_{2}a_{4}-a_{3}a_{5}+a_{4}a_{6}) 
\nonumber \\
e_{5} & = & -e_{7} = 64(a_{0}a_{1}-a_{1}a_{2}+a_{2}a_{3}-a_{3}a_{4}+a_{4}a_{5}-a_{5}a_{6}) \nonumber \\
e_{6} & = & 64(a_{0}^{2}-a_{1}^{2}+a_{2}^{2}-a_{3}^{2}+a_{4}^{2}-a_{5}^{2} + 
a_{6}^{2}). \nonumber 
\end{eqnarray}

\medskip 

\noindent
{\bf Step 3} computes the quotient of $E(x)$,  produced in 
Step 2, and $A(x)$, the denominator of the original integrand. In
the example discussed now we 
obtain
\begin{equation}
Z(x) = 64(a_{6}x^{6}-a_{5}x^{5}+a_{4}x^{4}-a_{3}x^{3}+a_{2}x^{2}-a_{1}x 
+ a_{0}). 
\end{equation}

\medskip

\noindent
{\bf Step 4} simply evaluates the product $C(x) = B(x)Z(x)$, where $B(x)$ is
the numerator of the original integrand and $Z(x)$ comes from Step 3. The 
polynomial $C(x)$ is written as 
\begin{equation}
C(x) = \sum_{k=0}^{10}c_{k}x^{10-k},
\end{equation}
\noindent
with
\begin{eqnarray}
c_{0} & = & 64a_{6}b_{0} \nonumber \\
c_{1} & = & -64(a_{5}b_{0}-a_{6}b_{1}) \nonumber \\
c_{2} & = & 64(a_{4}b_{0}-a_{5}b_{1}+a_{6}b_{2}) \nonumber \\
c_{3} & = & -64(a_{3}b_{0} - a_{4}b_{1} + a_{5}b_{2}-a_{6}b_{3}) \nonumber \\
c_{4} & = & 64(a_{2}b_{0}-a_{3}b_{1}+a_{4}b_{2}-a_{5}b_{3} + a_{6}b_{4}) 
\nonumber \\
c_{5} & = & -64(a_{1}b_{0}-a_{2}b_{1}+a_{3}b_{2}-a_{4}b_{3} + a_{5}b_{4}) 
\nonumber \\
c_{6} & = & 64(a_{0}b_{0}-a_{1}b_{1}+a_{2}b_{2}-a_{3}b_{3} + a_{4}b_{4}) 
\nonumber \\
c_{7} & = & 64(a_{0}b_{1} - a_{1}b_{2} + a_{2}b_{3}-a_{3}b_{4}) \nonumber \\
c_{8} & = & 64(a_{0}b_{2}-a_{1}b_{3}+a_{2}b_{4}) \nonumber \\
c_{9} & = & 64(a_{0}b_{3}-a_{1}b_{4}) \nonumber \\
c_{10} & = & 64a_{0}b_{4}. \nonumber 
\end{eqnarray}

\medskip

\noindent
{\bf Step 7} combines the functions defined in Step 5 and 6 to produce the 
new numerator 
\begin{equation}
J(x) = \sum_{k=0}^{4} d_{k}x^{4-k},
\end{equation}
\noindent
with
\begin{eqnarray}
d_{0} & = & 32(a_{6}b_{0}+a_{0}b_{4}) \nonumber \\
d_{1} & = & -16(a_{5}b_{0}-a_{6}b_{1}+a_{0}b_{3}-a_{1}b_{4}) \nonumber \\
d_{2} & = & 8(a_{4}b_{0}+3a_{6}b_{0}-a_{5}b_{1}+a_{0}b_{2}+a_{6}b_{2}-a_{1}b_{3}+3a_{0}b_{4}+a_{2}b_{4}) \nonumber \\
d_{3}& = & -4(a_{3}b_{0}+2a_{5}b_{0}+a_{0}b_{1}-a_{4}b_{1}-2a_{6}b_{1}-a_{1}b_{2}+a_{5}b_{2}) - \nonumber \\
 & & -4(2a_{0}b_{3}+a_{2}b_{3}-a_{6}b_{3}-2a_{1}b_{4}-a_{3}b_{4}) \nonumber \\
d_{4} & = & 2(a_{0}b_{0}+a_{2}b_{0}+a_{4}b_{0}+a_{6}b_{0}-a_{1}b_{1}-a_{3}b_{1}-a_{5}b_{1}) +  \nonumber \\
& & + 2(a_{0}b_{2}+a_{2}b_{2}+a_{4}b_{2}+a_{6}b_{2}-a_{1}b_{3}-a_{3}b_{3}+a_{0}b_{4}+a_{2}b_{4}+a_{4}b_{4}+a_{6}b_{4}). \nonumber 
\end{eqnarray}

\end{example}

\medskip

\begin{Note}
Given a rational function 
of order $p$ and a choice of method of order $m$, the calculation of $H$ 
and $J$ illustrated here is done \emph{once}.  We have produced a 
transformation sending 
\begin{equation}
F(x) = \frac{b_{0}x^{4} + b_{1}x^{3} + b_{2}x^{2} + b_{3}x + b_{4} }
{a_{0}x^{6} + a_{1}x^{5} + a_{2}x^{4} + a_{3}x^{3} 
+ a_{4}x^{2} + a_{5}x + a_{6} }
\end{equation}
\no
to
\begin{equation}
F_{1,6}(x) := \mathfrak{L}_{6,2}F(x) 
= \frac{d_{0}x^{4} + d_{1}x^{3} + d_{2}x^{2} + d_{3}x + d_{4} }
{h_{0}x^{6} + h_{1}x^{5} + h_{2}x^{4} + h_{3}x^{3} 
+ h_{4}x^{2} + h_{5}x + h_{6} }
\end{equation}
\noindent
with the new coefficients given as above. 
We expect to produce a \emph{precomputed} 
array of formulas, indexed by $(p,m)$, to be made available to the 
community. 
\end{Note}

In the next series of examples, we will assume that the formulas for the 
Landen transformations have been precomputed. 

\medskip

\begin{example}
The rational function 
\begin{equation}
F_{0}(x) = \frac{1}{x^{2} + 4x+15} 
\end{equation}
\no
satisfies 
\begin{equation}
I := \int_{-\infty}^{\infty} F_{0}(x) \, dx = \frac{\pi}{\sqrt{11}},
\end{equation}
\noindent
with numerical value 
\begin{equation}
I \sim  0.94722582509948293643. 
\end{equation}
\noindent
We now employ the algorithm, with a method of order $2$, to obtain the sequence
of rational functions 
\begin{equation}
F_{n,2}(x) = \frac{c_{n}}{x^{2} + a_{n}x + b_{n}},
\end{equation}
\noindent
with $F_{0,2}(x) = F_{0}(x)$, and the property that 
\begin{equation}
I = \int_{-\infty}^{\infty} F_{n,2}(x) \, dx = 
c_{n}\int_{-\infty}^{\infty} \frac{dx}{x^{2} + a_{n}x + b_{n}}.
\end{equation}
\noindent
The convergence analysis described in \cite{manna-moll2} shows that 
$a_{n} \to 0, \, b_{n} \to 1$, thus 
\begin{equation}
I = \lim\limits_{n \to \infty} \pi \, c_{n}.
\end{equation}
\noindent
Even though the limiting value of the integral depends only upon the terms 
$c_{n}$, the formulas to generate these values also involve
$a_{n}$ and $b_{n}$. Therefore one must store the current value of all 
the parameters. The
first few of them are shown in the next table: \\

\begin{table}[h]
\begin{center}
\begin{tabular}{|c|c|c|c|} \hline 
\hline 
$n$ & $c_{n}$ & $a_{n}$ & $b_{n}$   \\ \hline 
 \hline 
  & & &  \\
$0$ & $1$ & $4$ & $15$   \\ \hline 
  & & &  \\
$1$ & $\tfrac{8}{15}$ & $\tfrac{28}{15}$ & $4$ \\ 
& & & \\
\hline 
  & & & \\
$2$ & $\tfrac{1}{3}$ & $\tfrac{7}{10}$ & $\tfrac{4841}{3600}$ \\ 
  & & & \\
\hline 
  & & & \\
$3$ & $\tfrac{8441}{29046}$ & $\tfrac{8687}{96820}$ & 
$\tfrac{64900081}{69710400}$ \\ 
  & & &  \\
\hline
\end{tabular}
\end{center}
\vskip 0.2in
\caption{Rational Landen of order $2$}
\end{table}

\medskip

\begin{Note}
The expression $\pi c_{n}$ gives an approximation to the integral of $F_{0}(x)$.
For example, for $n=6$, we find that
\begin{equation}
 c_{6} = \frac{3471070386673821384824326347489289738211683509253931254760471}{11512238093504492278949475398059063785494372327433955614454608}
\nonumber
\end{equation}
\noindent
and $\pi c_{6}$ agrees with the integral of $F_{0}$ up to $13$ digits. 
The coefficients $(a_{n}, \, b_{n}, \, c_{n})$ are rational numbers and their
height increases with $n$. Recall that the \emph{height} 
of a rational number, written 
in irreducible form as $x = \frac{u}{v}$,  is 
\begin{equation}
h(x) = \text{Max} \{ |u|, \, |v| \}.
\end{equation}
\noindent
For example, at $n=10$, the heights of $a_{10}, \, b_{10}, \, c_{10}$ have 
approximately $1000$ digits. At this stage, the value $\pi c_{10}$ gives
$196$ digits of the integral. Naturally, this 
increases the complexity of the calculations
\emph{if we use exact arithmetic}. An interesting way out of this 
problem is to replace the rational number $c_{n}$ by the truncation of its 
continued fraction. For example, the first $20$ terms of the continued 
fraction of $c_{10}$, a rational number of height \emph{only} $6$, 
differs from $c_{10}$ by less than $10^{-24}$.  Details will be given in a
future publication. 
\end{Note}
\end{example}

\medskip

\begin{example}
Mathematica 6.0 shows that if
\begin{equation}
F(x) = \frac{1}{x^{4}+6x^{3}+16x^{2}+21x+13} 
\end{equation}
\no
then
\begin{equation}
I:= \int_{-\infty}^{\infty} F(x) \, dx =
2 \pi \,  \sqrt{ \frac{2}{111} ( \sqrt{37}-5) }.
\label{int-I}
\end{equation}
\no
We employ the identity (\ref{newint}) to map the problem to the 
interval $[0, \, 1]$.
The new rational function 
\begin{eqnarray}
g(x) & := & \frac{60x^{6} -288x^{5}+584x^{4}-648x^{3}+422x^{2}-156x+26}
{(3x^{4}-15x^{3}+31x^{2}-31x+13)(57x^4-153x^3+157x^2-73x+13)} \nonumber  \\
& = & \frac{60x^{6} -288x^{5}+584x^{4}-648x^{3}+422x^{2}-156x+26}
{171x^{8} -1314x^{7}+4533x^6-9084x^5+11485x^4-9314x^3+4707x^2-1352x+169}
\nonumber 
\end{eqnarray}
\noindent
satisfies 
\begin{equation}
\int_{0}^{1} g(x) \, dx = \int_{-\infty}^{\infty} F(x) \, dx.
\end{equation}
\noindent
We now compare the numerical approximation to the integral of $g$ over $[0,1]$ 
computed with the methods described here, with a numerical integration 
using the trapezoidal rule. A more systematic comparison with more 
sophisticated classical numerical schemes will be described elsewhere. 

The trapezoidal rule states that 
\begin{equation}
\int_{a}^{b} g(x) \, dx = \frac{h}{2}(g(a)+g(b)) - \frac{1}{12}h^{3} 
g''(\xi), 
\label{app}
\end{equation}
\noindent 
where $h = b-a$ and $\xi \in [a, \, b]$. Define
\begin{equation}
M := \text{Max} \{ |g''(t)|: \, a \leq t \leq b \}.
\end{equation}
\noindent
Then the \emph{error term}, namely the second term in (\ref{app}) is
bounded by $Mh^{3}/12$. To 
obtain an approximation to 
the integral $I$ in (\ref{int-I}), we choose 
$n \in \mathbb{N}$, partition $[0, \, 1]$ into $n$
intervals of equal  length $h = 1/n$ and apply the trapezoidal rule to each 
subinterval. This yields the expression 
\begin{equation}
I_{n} := \int_{0}^{1} g(x) \, dx \sim  \frac{1}{2n} \left( g_{0} + 
\sum_{i=1}^{n-2} g_{i} + g_{n} \right) - 
\frac{1}{12 n^{3}} \sum_{i=1}^{n-1} g''(\xi_{i}), 
\end{equation}
\noindent 
where $g_{i} = g(i/n)$ and  $\frac{i-1}{n} \leq \xi_{i} \leq \frac{i}{n}$. 
The total error in this approximation is bounded by $M/12n^{2}$. To compute 
the approximation $I_{n}$ to $I$ requieres the $n$ values $\{ g_{i}: \, 
0 \leq i \leq n \}$. The relative error $(I-I_{n})/I$ for $n=100$ is 
$5.29805 \times 10^{-6}$. It drops to $3.1505 \times 10^{-8}$ for $n=1000$ and 
to $2.9445 \times 10^{-10}$ for $n=10000$. 

\medskip

We now use the method of rational Landen transformations of order $m$
to produce 
approximations to the integral $I$ of $F(x)$ over $\mathbb{R}$.  Recall 
that the method yields a family of rational functions $R_{n,m}(x)$ with
integral $I$. For 
example, the first two functions for a method of order $2$ are
\begin{eqnarray}
R_{1,2}(x) & = & 
\frac{4(2x^{2}+6x+15)}{208x^{4}+456x^{3}+600x^{2} + 396x +171}, \nn \\
R_{2,2}(x) & = & \frac{8(13848x^{2}+11652x+11531)}{569088x^{4}-35136x^{3}+756384x^{2}-8616x+232537}. \nn
\end{eqnarray}
\noindent
The approximations to $I$ are then obtained from 
\begin{equation}
\text{app}_{n,m} := 
\frac{\text{ Constant term in the numerator of } R_{n,m}(x) } 
{\text{ Constant term in the denominator of } R_{n,m}(x) } \times \pi.
\end{equation}
\noindent
The next table shows the relative errors 
\begin{equation}
\text{rel}_{n,m} := \frac{|\text{app}_{n,m} - I|}{|I|}
\end{equation}
\noindent
for $ 2 \leq m \leq 6$ and $ 2 \leq n \leq 5$. 

\begin{table}[h]
\begin{center}
\begin{tabular}{|c|c|c|c|c|c|} \hline 
\hline 
  & & & & &  \\
$n$ & $m=2$ & $m=3$ & $m=4$ & $m=5$  & $m=6$ \\ 
  & & & & &  \\
 \hline \hline
  & & & & &  \\
$2$ & $0.30314$ & $0.022076$ & $0.0021170$ & $2.2646 \times 10^{-6}$  & $6.3257 \times 10^{-7}$  \\ \hline 
  & & &  & & \\
$3$ & $0.058475$ & $0.000035272$ & $5.2932 \times 10^{-12} $ 
& $2.9440 \times 10^{-23}$  & $4.4813 \times 10^{-40}$  \\ \hline 
  & & & & &  \\
$4$ & $0.0021170$ & $3.2713 \times 10^{-15}$ & $2.0616 \times 10^{-47} $ 
& $1.9758 \times 10^{-115}$  & $3.6655 \times 10^{-239}$  \\ \hline 
  & & & & &  \\
$5$ & $3.2700 \times 10^{-6}$ & $3.6952 \times 10^{-45}$ 
& $5.3750 \times 10^{-190} $ 
& $3.1671 \times 10^{-577}$  & $4.0442\times 10^{-1434}$  \\ \hline 
\hline 
\end{tabular}
\end{center}
\vskip 0.2in
\caption{Relative error for the numerical evaluation of $I_{n}$.}
\end{table}

\medskip

This table contains clear evidence to support the convergence orders claimed 
in (\ref{2point8}). 
\end{example}

\medskip

\begin{example}
The rational Landen transformations can be used to evaluate
\begin{equation}
\int_{-\infty}^{\infty} \frac{dx}{(x-2)^2 + \eps^2} = \frac{\pi}{\eps},
\label{int-eps}
\end{equation}
\noindent
for $\eps > 0$ small. This example illustrates the fact that the proposed 
method converges, even when the integrand has poles very close to the real 
axis. A systematic description of the sensitivity of the iteration with 
respect to the parameter $\eps$, will be presented 
elsewhere. 

For fixed $\eps > 0$, we apply a method of order $2$ to (\ref{int-eps}). 
This produces a sequence of rational functions of the form
\begin{equation}
R_{n}(x) = \frac{b_{0,n}}{a_{0,n} + a_{1,n}x + a_{2,n}x^{2}}
\end{equation}
\noindent
that satisfy 
\begin{equation}
\int_{-\infty}^{\infty} R_{n}(x) \, dx = \frac{\pi}{\eps}.
\end{equation}

The explicit Landen transformation of order $2$ is given by
\begin{eqnarray}
a_{0,n+1} & = & 
(a_{0,n}-a_{1,n}+a_{2,n})(a_{0,n}+a_{1,n}+a_{2,n}) \label{landen-1} \\
a_{1,n+1} & = & 2a_{1,n}(a_{0,n}-a_{2,n}) \nonumber \\
a_{2,n+1} & = & 4a_{0,n}a_{2,n} \nonumber \\
b_{0,n+1} & = & 2b_{0,n}(a_{0,n}+a_{2,n}), \nonumber
\end{eqnarray}
\noindent
with initial conditions 
\begin{equation}
a_{0,0} = 4 + \eps^{2}, \, a_{1,0} = -4, \, a_{2,0} = 1, \, b_{0,0} = 1.
\end{equation}
\noindent 

The sequence $R_{n}(x)$ has coefficients that depend upon the parameter 
$\eps$. For example, 
\begin{equation}
R_{1}(x) = \tfrac{2(5+ \eps^{2})}{4(4 + \eps^{2})x^{2} -8(3 + \eps^{2})x + 
(1+\eps^{2})(9+\eps^{2})} 
\nonumber 
\end{equation}
\noindent
and 
\begin{equation}
R_{2}(x) = \tfrac{4(5+\eps^{2})(25 + 14 \eps^{2} + \eps^{4})}
{(1+ 6 \eps^{2} + \eps^{4})(49 + 22 \eps^{2} + \eps^{4}) -16(-1+\eps)
(1 + \eps)(3 + \eps^2)(7 + \eps^2)x + 16 (1 + \eps^{2})(4 + \eps^{2})
(9 + \eps^2)}. 
\nonumber
\end{equation}

The theory described above shows that, for fixed $\eps > 0$ and $n \to \infty$, the sequences 
\begin{equation}
\frac{a_{0,n}}{b_{0,n}} \to L, \quad  
\frac{a_{1,n}}{b_{0,n}} \to 0, \quad 
\frac{a_{2,n}}{b_{0,n}} \to L 
\end{equation}
\noindent
converge to the stated limits. Moreover, the 
invariance of (\ref{int-eps}) under the transformations given 
in (\ref{landen-1}) show that $L = \eps$. 

Define the error 
\begin{equation}
\text{err}_{n} := \left( \left( \tfrac{a_{0,n}}{b_{0,n}} - \eps \right)^{2} +
\left( \tfrac{a_{1,n}}{b_{0,n}} \right)^{2} +
\left( \tfrac{a_{2,n}}{b_{0,n}} - \eps \right)^{2} \right)^{1/2}.
\end{equation}
\noindent
Then Table 3 shows the ratios $\text{err}_{16}/\text{err}_{15}$ obtained 
after $15$ iterations of (\ref{landen-1}) for methods
of order $2$ and $3$. The 
calculations are done with $10^6$ digit precision.

\begin{table}[ht]
\begin{center}
\begin{tabular}{|c|c|c|} \hline 
\hline 
  &  &  \\
$\eps$ & $\text{order } 2$  & $\text{order } 3$\\ 
  &  &  \\
 \hline \hline
  &  &  \\
$.1$ & $3.58047 \times 10^{-569}$ & $3.49118 \times 10^{-497802}$  \\ \hline 
  & &  \\
$.01$ & $1.36862 \times 10^{-57}$ & $4.24935 \times 10^{-49853}$  \\ \hline 
  & &  \\
$.001$ & $2.07254 \times 10^{-6} $ & $4.73905 \times 10^{-4986}$  \\ \hline 
  & &  \\
$.0001$ & $2.16805 \times 10^{-2}$ & $3.48094 \times 10^{-499}$  \\ \hline 
  & &  \\
$.00001$ & $4.68150 \times 10^{-1}$ & $1.62880 \times 10^{-50}$  \\ \hline 
\end{tabular}
\end{center}
\vskip 0.2in
\caption{The quotient $\text{err}_{16}/\text{err}_{15}$ as a function of 
the parameter $\eps$.}
\end{table}

\medskip

The data in Table $3$ shows the exponential decay of the error. 
Given a tolerance $\delta > 0$,  we have 
observed that the number of steps required to achieve $\text{err}_{n} < 
\delta$ increases as $\eps \to 0$. A quantitative description of this 
phenomena is in preparation and it will be reported elsewhere. 

\end{example}

\medskip

\section{Conclusions}

We have described the rational Landen transformations and their use in 
the numerical integration of rational functions.  We have exhibited fast 
convergence of this method and presented an example comparing it to 
the classical integration schemes. 

A systematic comparative analysis of this method with respect to standard 
numerical algorithms will be discussed elsewhere. An interesting challenging 
problem is to extend the use of rational Landen transformations to produce 
fast numerical integrators for arbitrary functions. In particular the method 
is well suited for the numerical integration of meromorphic function with 
poles off the real line. 

\medskip

\no
{\bf Acknowledgments}. The first author is partially funded by the 
AARMS Director's Postdoctoral Fellowship.  The work of the second 
author was partially funded by
$\text{NSF-DMS } 0070567$. The authors wish to thank a referee for suggesting
the integral described in the last example.

\end{document}